\documentclass{birkjour_t2}
 \theoremstyle{definition}
 
 \theoremstyle{remark}

 \numberwithin{equation}{section}

\usepackage{graphicx}
\usepackage{amssymb, amsmath, amsthm}
\usepackage{amsfonts}
\usepackage{amssymb}
\usepackage{float}
\usepackage{mathrsfs}
\usepackage{enumitem} 
\newcommand{\dis}{\displaystyle}

\newcommand{\divv}{\text{\rm div}}

\newcommand{\R}{\mathbb R}

\renewcommand{\theequation}{\thesection.\arabic{equation}}
\numberwithin{equation}{section}
\numberwithin{theorem}{section}

\numberwithin{figure}{section}

\begin{document}

%
%
%
%
%
%
%
%
%

\title[Navier-Sokes Equations with Large External Force]
 {Global Solutions of the Navier-Stokes Equations\\
 for Isentropic Flow with Large External Potential Force}

\author[Anthony Suen]{Anthony Suen}

\address{%
Department of Mathematics\\Indiana
    University\\Bloomington, IN 47405}

\email{cksuen@indiana.edu}



\date{September 27, 2011}
\dedicatory{To my family and my wife Candy}

\begin{abstract}
We prove the global-in-time existence of weak solutions to the Navier-Stokes equations of compressible isentropic flow in three space dimensions with adiabatic exponent $\gamma\ge1$. Initial data and solutions are small in $L^2$ around a non-constant steady state with densities being positive and essentially bounded. No smallness assumption is imposed on the external forces when $\gamma=1$. A great deal of information about partial regularity and large-time behavior is obtained.
\end{abstract}

\maketitle
\section{Introduction}
We prove the global existence of weak solutions to the Navier-Stokes equations of compressible flow in three space dimensions:
\begin{align}
\label{1.1} \left\{ \begin{array}{l}
\rho_t + \divv (\rho u) =0, \\
(\rho u^j)_t + \divv (\rho u^j u) + (P_{\gamma})_{x_j} = \mu\,\Delta u +
(\xi - \mu) \, (\divv \,u)_{x_j}  + f^{j}.
\end{array}\right.
\end{align}
Here $\rho$ and $u=(u^1,u^2,u^3)$ are the unknown functions of $x\in\R^3$ and $t\ge0$, $P_{\gamma}=P_{\gamma}(\rho)$
is the
pressure, $f$ is the external force, $\mu$ and $\xi-\mu$ are viscosity constants.

The system \eqref{1.1} is solved subject to initial conditions
\begin{equation}
\label{1.2}
(\rho(\cdot,0),u(\cdot,0)) = (\rho_0,u_0),
\end{equation}
where $\rho_0$ is bounded above and below away from zero,  $u_0\in L^p$ for some   $p>6$, and modulo constants, $(\rho_0-\tilde\rho_{\gamma},u_0)$ is small in $L^2(\R^3)$ for some specific non-constant function $\tilde\rho_{\gamma}$ which will be defined later. The solvability to various Cauchy problems for the Naiver-Stokes equations has been discussed by many other mathematicians for decades. Matsumura-Nishida \cite{mn1} proved the global existence of $H^3$ solutions around a constant state when the initial data is taken to be small in $H^3$, and later Danchin \cite{danchin} generalized their results by replacing $H^3$ with certain Besov spaces of functions. On the other hand, Lions \cite{lions} and Feireisl \cite{feireisl1}-\cite{feireisl2} proved global existence of weak solutions to \eqref{1.1}-\eqref{1.2} with finite-energy initial data and nonnegative initial density. In between those two type of solutions as mentioned, Hoff \cite{hoff1}-\cite{hoff4} developed a new theory of {\em intermediate regularity} class solutions which may exhibit discontinuities in density and velocity gradient across hypersurfaces in $\R^3$. In the presence of large external force, our work generalizes earlier results of Matsumura-Yamagata \cite{my} in two ways: the restriction on the $L^\infty$ norm of $\rho-\tilde\rho_\gamma$ has been eliminated and initial velocity is not necessary in $H^1$.

\medskip

We introduce two variables associated with the system \eqref{1.1} which are important to our analysis. The first one is the usual vorticity matrix $\omega=\omega^{j,k}=u^j_{x_k}-u^k_{x_j}$, while the other one is the effective viscous flux $F$ given by 
\begin{equation}\label{1.3}
F=\xi\divv\,u-P_\gamma(\rho)+P_\gamma(\tilde\rho_\gamma).
\end{equation}
By adding and subtracting terms, we can rewrite the momentum equation as in \eqref{1.1} in terms of $F$ and $\omega$:
\begin{equation}\label{1.4}
\rho\dot u^j=F_{x_j}+\mu\omega^{j,k}_{x_k}+\rho f^j-P_\gamma(\tilde\rho_\gamma)_{x_j}.
\end{equation}
The decomposition \eqref{1.4} also implies that
\begin{equation}\label{1.5}
\Delta F=\divv(\rho\dot u+\rho f-\nabla P_\gamma(\tilde\rho_\gamma)).
\end{equation}
We refer to Hoff \cite{hoff1} for a more detailed discussion of $F$. 
\medskip

We now give a precise formulation of our results. First concerning the pressure $P_\gamma$ we assume that

\begin{enumerate}[label={\upshape{(1.\arabic*)}}, ref={\upshape{2.}\arabic*}, topsep=0.3cm, itemsep=0.15cm] 
\setcounter{enumi}{5}
\item there is $L>0$ such that $P_{\gamma}(\rho)=L\rho^{\gamma}$ for $\gamma\ge1$.
\end{enumerate}
\setcounter{equation}{6}
We also fix a positive reference density $\rho_\infty$ and choose bounding densities $0<\underline\rho<\bar\rho$ satisfying
\begin{equation}\label{1.6}
\underline\rho<\rho_\infty<\bar\rho
\end{equation}
and we define a positive number $\delta$ such that
\begin{equation}\label{1.7}
\delta=\frac{1}{2}\min\{|\rho_\infty-\underline\rho|,|\rho_\infty-\bar\rho|\}.
\end{equation}

Concerning the external force $f$, we assume that
\begin{enumerate}[label={\upshape{(1.\arabic*)}}, ref={\upshape{2.}\arabic*}, topsep=0.3cm, itemsep=0.15cm]
\setcounter{enumi}{8} 
\item $f=-\rho\nabla\psi$, where $\psi\in H^4(\R^3)$;
\item there is a constant $\mathcal C>0$ such that $|\psi(x)|+|\nabla\psi(x)|\le\mathcal C$ and $|D_{x}^2\psi(x)|\le|\nabla\psi(x)|$ for $x\in\R^3$;
\item $\lim\limits_{|x|\rightarrow\infty}\psi(x) = 0$.
\end{enumerate}
\setcounter{equation}{11}

Concerning the diffusion coefficients $\mu$ and $\xi$, we assume that
\begin{equation}\label{1.1.15}
0<\xi<\Big(\textstyle\frac{3}{2}+{\textstyle{\frac{\sqrt{21}}{6}\Big)\mu}}.
\end{equation}
It follows that
\begin{equation}\label{1.1.16}
{\textstyle\frac{1}{4}}\mu(p-2)-\frac{[{\textstyle\frac{1}{4}}(\xi-\mu)(p-2)]^2}{\frac{1}{3}\mu+(\xi-\mu)}>0
\end{equation}
for $p=6$ and consequently for some $p>6$, which we now fix.
\medskip

We define $\tilde\rho_\gamma$ as mentioned at the beginning of this section. Given a positive constant densty $\rho_{\infty}$, we say $(\tilde\rho_\gamma,0)$ is a {\it steady state solution} to \eqref{1.1} if $\tilde\rho_\gamma\in C^2(\R^3)$ and the following holds
\begin{align}
\label{1.1.17} \left\{ \begin{array}{l}
\nabla P_\gamma(\tilde\rho(x)) =-\tilde\rho_\gamma(x)\nabla\psi(x), \\
\lim\limits_{|x|\rightarrow\infty}\tilde\rho(x) = \rho_{\infty}.
\end{array}\right.
\end{align}
By direct computation, $\rho_\gamma$ has the following explicit form
$$\rho_\gamma(x)=\left\{
\begin{array}{ll}
\rho_\infty\exp\left[-\frac{1}{L}\psi(x)\right]                     &\mbox{if $\gamma=1,$} \\
\left[\rho_\infty^{\gamma-1}-\frac{\gamma-1}{L\gamma}\psi(x)\right]^{\frac{1}{\gamma-1}}                     &\mbox{if $\gamma>1.$} \\
\end{array}
\right.\leqno{(1.15)}$$
\medskip

Concerning the initial data $(\rho_0,u_0)$ we assume that there is a positive number $N$, which may be arbitrarily large such that
\setcounter{equation}{15}
\begin{equation}\label{1.2.1}
\|u_0\|_{L^{p}}\le N,
\end{equation}
and there is positive number $\delta$ with $d<\delta$ such that
\begin{equation}\label{1.2.2}
\underline\rho+d < {\rm ess}\inf\rho_0 \le {\rm ess}\sup\rho_0<\bar\rho -d,
\end{equation}
\rm
We also write
\hfill
\begin{equation}\label{1.2.35}
C_{0}=\|\rho_0-\tilde{\rho_{\gamma}}\|_{L^2}+\|u_0\|_{L^{2}}.
\end{equation}
where $\tilde\rho_\gamma$ is defined as in (1.15).
\medskip

Weak solutions are defined in the usual way; we say that $(\rho,u)$ is a weak solution of \eqref{1.1}-\eqref{1.2} provided that \nobreak{$(\rho-\tilde{\rho_\gamma},\,\rho u)\in C([0,\infty);H^{-1}(\R^3))$ with $(\rho,u)|_{t=0}=(\rho_0,u_0)$, $\nabla u\in L^2(\R^3\times(0,\infty))$ for $t>0$, and 
the following identities hold for times $t_2\ge t_1 \ge 0$ and $C^1$ test functions $\varphi$ having uniformly bounded support in $x$ for $t\in[t_1,t_2]$:
\begin{align}\label{1.2.5}
\left.\int_{\R^3}\rho(x,\cdot)\varphi(x,\cdot)dx\right|_{t_1}^{t_2}=\int_{t_1}^{t_2}\!\!\!\int_{\R^3}(\rho\varphi_t + \rho u\cdot\nabla\varphi)dxdt,
\end{align}
and
\begin{align}\label{1.2.6}
\left.\int_{\R^3}(\rho u^{j})(x,\cdot)\varphi(x,\cdot)dx\right|_{t_1}^{t_2}=\int_{t_1}^{t_2}\!\!\!&\int_{\R^3}[\rho u^{j}\varphi_t + \rho u^{j}u\cdot\nabla\varphi + P(\rho)\varphi_{x_j}]dxdt\notag\\
& - \int_{t_1}^{t_2}\!\!\!\int_{\R^3}[\mu\nabla u^{j}\cdot\nabla\varphi + (\mu - \xi)(\divv\,u)\varphi_{x_j}]dxdt\\&\qquad\qquad\qquad\qquad\qquad\,\,\,+ \int_{t_1}^{t_2}\!\!\!\int_{\R^3}\rho\varphi_{x_j}\psi dxdt.\notag
\end{align}

We use the usual notation for H\"older seminorms: for $v:\R^3\to \R^m$ and $\alpha \in (0,1]$, 
\hfill
$$\langle v\rangle^\alpha = \sup_{{x_1,x_2\in 
\R^3}\atop{x_1\not=x_2}}
{{|v(x_2) -v(x_1)|}\over{|x_2-x_1|^\alpha}}\,;$$
and for $v:Q\subseteq\R^3 \times[0,\infty)\to \R^m$ and $\alpha_1,\alpha_2 \in (0,1]$,
\hfill
$$\langle v\rangle^{\alpha_1,\alpha_2}_{Q} = \sup_{{(x_1,t_1),(x_2,t_2)\in 
Q}\atop{(x_1,t_1)\not=(x_2,t_2)}}
{{|v(x_2,t_2) - v(x_1,t_1)|}\over{|x_2-x_1|^{\alpha_1} + |t_2-t_1|^{\alpha_2}}}\,.$$
We denote the material derivative of a given function $v$ by
$\dot{v}=v_t + \nabla v\cdot u$,
and if $X$ is a Banach space we will abbreviate $X^3$ by $X$ when convenient. If $I\subset [0,\infty)$ is an interval, $C^1(I;X)$ will be the elements $v\in C(I;X)$ such that the distribution derivative $v_t\in {\mathcal D}'(\R^3\times{\rm int}\,I)$ is realized as an element of $C(I;X)$. Finally, if $\Omega\subset\R^3$ is a Lebesgue measurable subset in $\R^3$, $|\Omega|$ will be the corresponding volume of $\Omega$.

\medskip
 
The following is the main result of this paper, which is valid for the case when $\gamma=1$:

\medskip

\noindent{\bf Theorem~1.1} \em Let $\gamma=1$ and assume that the system parameters in {\rm \eqref{1.1}} satisfy the conditions in {\rm (1.6) and (1.9)-\eqref{1.1.16}}. Let positive numbers $\mathcal{C}, L,\rho_{\infty},\bar\rho,\underline\rho,\delta$ be given satisfying \eqref{1.6}-\eqref{1.7} and 
\begin{equation}\label{1.9.1}
L\log\left(\frac{\rho_\infty}{\bar\rho}\right)<-\mathcal C<\mathcal C<L\log\left(\frac{\rho_\infty}{\underline\rho}\right)
\end{equation}
Then given positive numbers $N$ and $d<\delta$, there are positive constants $a, C, \theta$ depending on the parameters and assumptions in {\rm (1.6) and (1.9)-\eqref{1.1.16}}, on $\mathcal{C}, L,\rho_{\infty},\bar\rho,\underline\rho,R$, on $N$ and on a positive lower bound for $d$, such that if an initial data $(\rho_0,u_0)$ is given satisfying \eqref{1.2.1}-\eqref{1.2.35} with
\begin{equation}\label{1.3.1}
C_0<a,
\end{equation}
then there is a solution $(\rho,u)$ to \eqref{1.1}-\eqref{1.2} in the sense of \eqref{1.2.5}-\eqref{1.2.6} on all of $\R^3\times[0,\infty)$. The solution satisfies the following:
\begin{equation}\label{1.4.1}
\rho-\tilde{\rho}_\gamma,\,\rho u\in C([0,\infty);H^{-1}(\R^3)),
\end{equation}
\begin{equation}\label{1.4.2}
\nabla u\in L^2(\R^3\times(0,\infty)),
\end{equation}
\begin{equation}\label{1.4.3}
u(\cdot,t)\in H^1 (\R^3),\;t>0,
\end{equation}
\begin{equation}\label{1.4.4}
\mbox{$\omega(\cdot,t),F(\cdot,t)\in H^1 (\R^3)$, $t>0$,}
\end{equation}
\begin{equation}\label{1.4.6}
\langle u\rangle^{\frac{1}{2},\frac{1}{4}}_{\R^3 \times [\tau,\infty)}\leq C(\tau)C_{0}^{\theta},\;\tau>0,
\end{equation}
where $C(\tau)$ may depend additionally on a positive lower bound for $\tau$,
\begin{equation}\label{1.4.7}
\mbox{$\underline\rho\le\rho(x,t)\le\bar\rho$ a.e. on $\R^3\times[0,\infty)$},
\end{equation}
and
\begin{align}\label{1.4.8}
\sup_{t>0}\int_{\R^3}\big[(\rho - \tilde\rho_\gamma)^2 + |u|^2 &+ \sigma|\nabla u|^2 + \sigma^3 ( F^2 + |\nabla\omega|^2 )\big]dx\notag\\
+\int_{0}^{\infty}\!\!\!\int_{\R^3}\big[|\nabla u|^2 +& \sigma(|\dot u|^2 + |\nabla\omega|^2 )+\sigma^{3}|\nabla\dot{u}|^2\big]dxds\le CC_{0}^{\theta}&
\end{align}
where $\sigma(t)=\min\{1,t\}$. Moreover, $(\rho,u)\to(\tilde\rho_\gamma,0)$ as $t\to\infty$ in the sense that, for all $r_1\in(2,\infty]$ and $r_2\in(2,\infty)$,
\begin{align}\label{1.30}
\lim\limits_{t\rightarrow\infty}||\rho(\cdot,t)-\tilde\rho_\gamma(\cdot)||_{L^{r_1}}+||u(\cdot,t)||_{L^{r_2}}=0.
\end{align}
\rm
\bigskip

Using similar method, we also obtain parallel results for $\gamma>1$, with an extra assumption on the support of the external potential force $\psi$:

\medskip

\noindent{\bf Theorem~1.2} \em Let $\gamma>1$ and assume that the system parameters in {\rm \eqref{1.1}} satisfy the conditions in {\rm (1.6) and (1.9)-\eqref{1.1.16}}. Let positive numbers $\mathcal{C}, L,\rho_{\infty},\bar\rho,\underline\rho,\delta$ be given satisfying \eqref{1.6}-\eqref{1.7} and 
\begin{equation}\label{1.9.2}
\frac{L\gamma}{\gamma-1}\left[\rho_\infty^{\gamma-1}-\bar\rho^{\gamma-1}\right]<-\mathcal C<\mathcal C<\frac{L\gamma}{\gamma-1}\left[\rho_\infty^{\gamma-1}-\underline\rho^{\gamma-1}\right]
\end{equation}
Then given positive numbers $N$ and $d<\delta$, there are positive constants $a, C, \theta$ depending on the parameters and assumptions in {\rm (1.6) and (1.9)-\eqref{1.1.16}}, on $\gamma,\mathcal{C}, L,\rho_{\infty},\bar\rho,\underline\rho$, on $N$ and on a positive lower bound for $d$, such that if an initial data $(\rho_0,u_0)$ is given satisfying \eqref{1.2.1}-\eqref{1.2.35} and the small-energy assumption \eqref{1.3.1}, and if
\begin{equation}
|{\rm supp}(\psi)|<a,
\end{equation}
then there is a solution $(\rho,u)$ to \eqref{1.1}-\eqref{1.2} in the sense of \eqref{1.2.5}-\eqref{1.2.6} on all of $\R^3\times[0,\infty)$. The solution satisfies \eqref{1.4.1}-\eqref{1.4.6} with bounds \eqref{1.4.7}-\eqref{1.4.8}, and $(\rho,u)\to(\tilde\rho_\gamma,0)$ as $t\to\infty$ in the sense of \eqref{1.30} for all $r_1\in(2,\infty]$ and $r_2\in(2,\infty)$.
\rm
\bigskip

We point out that for the case when $\gamma=1$, system \eqref{1.1} becomes simpler than those for $\gamma>1$. For example, in deriving {\em a priori} bounds for the term $\dis\int_{0}^{t}\!\!\!\int_{\R^3}\sigma|\dot u|^2 dxds$, we first rewrite the momentum equation in \eqref{1.1} as follows
\begin{align}\label{1.99}
\rho\dot u&+\tilde P_\gamma\nabla(P_\gamma \tilde P_\gamma^{-1}-1)+\rho \tilde P_\gamma^{-1}\nabla \tilde P_\gamma(\rho^{-1}P_\gamma-\tilde\rho_\gamma^{-1}\tilde P_\gamma)-\mu\Delta u-(\xi-\mu)\nabla(\divv\,u)=0.
\end{align}
When $\gamma=1$, the third of the above reads
$$\rho \tilde P_\gamma^{-1}\nabla \tilde P_\gamma(\rho^{-1}P_\gamma-\tilde\rho_\gamma^{-1}\tilde P_\gamma)=\rho \tilde P_\gamma^{-1}\nabla \tilde P_\gamma(L-L)=0,$$
so if we multiply \eqref{1.99} with $\sigma\dot u$ and integrate, $\dis\int_{0}^{t}\!\!\!\int_{\R^3}\sigma|\dot u|^2 dxds$ can be bounded in terms of the $L^2$ norms of $(\rho-\tilde\rho,u)$ and the $H^1$ norms of $u$ (modulo higher order terms in $u$). 

On the other hand, when $\gamma>1$, without any cancellation, we then have to estimate the term $\dis\int_{0}^t\!\!\!\int_{\R^3}\sigma\dot u\rho \tilde P_\gamma^{-1}\nabla \tilde P_\gamma(\rho^{-1}P_\gamma-\tilde\rho_\gamma^{-1}\tilde P_\gamma)dxds$ which makes the analysis much more intricate and forces extra assumptions imposed on the external force. We will explain more later in subsequent sections.
\medskip

This paper is organized as follows. We begin the proofs of Theorem~1.1 and 1.2 in section~2 with a number of {\em a priori} bounds for local-in-time smooth solutions. Since many of these estimates are rather long and technical, we omit those which are identical to or nearly identical to arguments given elsewhere in the literature. In section~3 we derive the necessary bounds for density by applying the estimates in Theorem~3.1 and 3.2 in a maximum principle argument along particle trajectories of the velocity, making important use of the monotonicity of $P_\gamma$ as described in \rm{(1.8)}. The small-energy assumption \eqref{1.3.1} then enables us to close these arguments to show in Theorem 3.1 and 3.2 that both the pointwise bounds for density and the {\em a priori} energy bounds of Theorem~3.1 do hold as long as the smooth solution exists. Finally in section~4 we prove Theorem~1.1 and 1.2 by constructing weak solutions as limits of smooth solutions corresponding to mollified initial data.

\medskip

We make use of the following standard facts (see Ziemer \cite{ziemer} Theorem~2.1.4, Remark~2.4.3, and Theorem~2.4.4, also Ladyzhenskaya \cite{ladyzhenskaya} section~1.4, for example). First, given $r\in[2,6]$ there is a constant $C(r)$ such that for $w\in H^1 (\R^3)$,
\begin{equation}\label{1.5.1}
\|w\|_{L^r(\R^3)} \le C(r) \left(\|w\|_{L^2(\R^3)}^{(6-r)/2r}\|\nabla w\|_{L^2(\R^3)}^{(3r-6)/2r}\right).
\end{equation}
Next, for any $r\in (1,\infty)$ there is a constant $C(r)$ such that for $w\in W^{1,r}(\R^3)$,
\begin{equation}\label{1.5.2}
\|w\|_{L^\infty (\R^3)} \le C(r) \|w\|_{W^{1,r}(\R^3)}
\end{equation}
and
\begin{equation}\label{1.5.3}
\langle w\rangle^\alpha_{\R^3}\le C(r)\|\nabla w\|_{L^r(\R^3)},
\end{equation}
where $\alpha=1-3/r$. If $\Gamma$ is the fundamental solution of the Laplace operator on $\R^3$, then there is a constant $C$ such that for any $f\in L^2(\R^3)\cap L^4(\R^3)$, 
\begin{equation}\label{1.5.4}
\|\Gamma_{x_j}*f\|_{L^\infty (\R^3)} \le C\left(\|f\|_{L^2(\R^3)} + \|f\|_{L^4(\R^3)}\right).
\end{equation}
Finally, there is a constant $M$ such that for any $v\in H^1(\R^3)$,
\begin{equation}\label{1.5.5}
\int_{\R^3}\frac{|v(x)|^2}{(1+|x|)^2}dx\le M\int_{\R^3}|\nabla v|^2dx.
\end{equation}
\medskip

\section{Energy Estimates}
\bigskip

In this section we derive {\em a priori} bounds for smooth, local-in-time solutions of \eqref{1.1}-\eqref{1.2} whose densities are strictly posititve and bounded. Specifically, we fix a smooth solution $(\rho-\tilde\rho_\gamma,u)$ of \eqref{1.1}-\eqref{1.2} on $\R^3\times[0,T]$ for some time $T>0$ with initial data $(\rho_0-\tilde\rho_\gamma,u_0)$. These bounds will depend only on the quantities $C_0,\mathcal{C}, L,\rho_{\infty},\bar\rho,\underline\rho,N,d$ and will be independent of the initial regularity and the time of existence.  

We define a functional $A(t)$ for a given such solution by
\begin{align}\label{3.1}A(t)=\sup_{0<s\le t}\int_{\R^3}\left[\sigma|\nabla u|^2 + \sigma^3 (|\dot{u}|^2 + |\nabla\omega|^2)\right]dx\\ 
+\int_{0}^{t}\!\!\!\int_{\R^3}\left[\sigma(|\dot{u}|^2 + |\nabla\omega|^2) + \sigma^3|\nabla\dot{u}|^2\right]dxds,\notag\end{align}
where $\sigma(t)\equiv\min\{1,t\}$, and we obtain the following {\em a priori} bound for $A(t)$ under the assumptions that the initial energy $C_0$ in \eqref{1.2.35} is small and that the density remains bounded above and below away from zero when $\gamma=1$:

\medskip
\noindent{\bf Theorem~2.1} \em  Let $\gamma=1$. Assume that the system parameters in {\rm \eqref{1.1}} satisfy the conditions in {\rm (1.6) and (1.9)-\eqref{1.1.16}}. Let positive numbers $\mathcal{C}, L,\rho_{\infty},\bar\rho,\underline\rho,\delta$ be given satisfying \eqref{1.6}-\eqref{1.7} and \eqref{1.9.1}. Then given positive numbers $N$ and $d<\delta$, there are positive constants $a, M, \theta$ depending on the parameters and assumptions in {\rm (1.6) and (1.9)-\eqref{1.1.16}}, on $\mathcal{C}, L,\rho_{\infty},\bar\rho,\underline\rho$, on $N$ and on a positive lower bound for $d$, such that if $(\rho-\tilde\rho_\gamma,u)$ is a solution of \eqref{1.1}-\eqref{1.2} on $\R^3\times[0,T]$ with initial data $(\rho_0-\tilde\rho_\gamma,u_0)\in H^3(\R^3))$ satisfying \eqref{1.2.1}-\eqref{1.2.35} with $C_0 <a$, and if 
\begin{align*}
\mbox{$\underline\rho\le\rho(x,t)\le\bar\rho$ on $\R^3\times[0,T]$},
\end{align*}
then 
\begin{align*}
A(T)\le MC_{0}^{\theta}.
\end{align*}
\rm
\bigskip

With an extra assumption on supp$(\psi)$, we can also obtain the following estimate for $\gamma>1$:

\medskip

\noindent{\bf Theorem~2.2} \em  Let $\gamma>1$. Assume that the system parameters in {\rm \eqref{1.1}} satisfy the conditions in {\rm (1.6) and (1.9)-\eqref{1.1.16}}. Let positive numbers $\mathcal{C}, L,\rho_{\infty},\bar\rho,\underline\rho,\delta$ be given satisfying \eqref{1.6}-\eqref{1.7} and \eqref{1.9.2}. Then given positive numbers $N$ and $d<\delta$, there are positive constants $a, M, \theta$ depending on the parameters and assumptions in {\rm (1.6) and (1.9)-\eqref{1.1.16}}, on $\gamma,\mathcal{C}, L,\rho_{\infty},\bar\rho,\underline\rho$, on $N$ and on a positive lower bound for $d$ such that if $|{\rm supp}(\psi)|<a$ and if $(\rho-\tilde\rho_\gamma,u)$ is a solution of \eqref{1.1}-\eqref{1.2} on $\R^3\times[0,T]$ with initial data $(\rho_0-\tilde\rho_\gamma,u_0)\in H^3(\R^3))$ satisfying \eqref{1.2.1}-\eqref{1.2.35} with $C_0 <a$ and
\begin{align*}
\mbox{$\underline\rho\le\rho(x,t)\le\bar\rho$ on $\R^3\times[0,T]$},
\end{align*}
then 
\begin{align*}
A(T)\le MC_{0}^{\theta}.
\end{align*}
\rm
\medskip

The proof will be given  in a sequence of lemmas in which we estimate a number of auxiliary functionals. To describe these we first recall the definition \eqref{1.1.16} of $p$, which is an open condition, and which therefore allows us to choose $q\in(6,p)$ which also satisfies \eqref{1.1.16}.
Then for a given $(\rho-\tilde\rho_\gamma,u)$, we define
\begin{align}\label{3.2}
H(t)=\int_{0}^{t}\!\!\!\int_{\R^3}\Big[\sigma^{3/2}|\nabla u|^3 & + \sigma^3 |\nabla u|^4 \Big]dxds\notag
\\&+ \Big |\sum_{1\le k_i,j_m\le 3}\int_{0}^{t}\!\!\!\int_{\R^3}\sigma u^{j_1}_{x_{k_1}}u^{ j_2}_{x_{k_2}}u^{j_3}_{x_{k_3}}dxds\Big |,
\end{align}
\begin{align}\label{3.3}
D(t)=\int_{0}^{t}\!\!\!\int_{\R^3}|\rho-\tilde\rho_\gamma|^2|\nabla P(\tilde\rho)|^2dxds.
\end{align}
Also, we can readily see that assumptions \eqref{1.9.1} and \eqref{1.9.2} imply $\tilde\rho_\gamma$ is well-defined and
\begin{align}\label{1.9.3}
\underline\rho<\tilde\rho_{\gamma}<\bar\rho,
\end{align}
which will be crucial to the later analysis. 

For simplicity, we write $\sigma=\sigma(t)=\min\{1,t\}$, $P_\gamma=P_\gamma(\rho)$ and $\tilde P_\gamma=P_\gamma(\tilde\rho_\gamma)$ without further referring.

\medskip

We begin with the following $L^2$ energy estimate, which is valid for $\gamma\ge1$:
\bigskip

\noindent{\bf Lemma 2.3} \em Assume that the system parameters in {\rm \eqref{1.1}} satisfy the conditions in {\rm (1.6) and (1.9)-\eqref{1.1.16}}. Then if $(\rho-\tilde\rho_\gamma,u)$ is a solution of \eqref{1.1}-\eqref{1.2} on $\R^3\times[0,T]$ with initial data $(\rho_0-\tilde\rho_\gamma,u_0)\in H^3(\R^3))$ satisfying \eqref{1.2.1}-\eqref{1.2.35}, and if $\rho,\tilde\rho\in[\underline\rho,\bar\rho]$, then for $\gamma\ge1$,
\begin{equation}\label{3.4}
\sup_{0\le t\le T}\int_{\R^3}(|\rho - \tilde\rho_\gamma|^2 + \rho|u|^2)dx + \int_{0}^{T}\!\!\!\int_{\R^3}|\nabla u|^2dxdt\leq MC_{0}.
\end{equation}
\rm
\begin{proof}
We use $\nabla\psi=-\tilde\rho_\gamma^{-1}\nabla\tilde P_\gamma$ on the momentum equation to get
\begin{equation}\label{3.5}
\rho\dot u+\rho(\rho^{-1}\nabla P_\gamma-\tilde\rho_\gamma^{-1}\nabla \tilde P_\gamma)-\mu\Delta u-(\xi-\mu)\nabla(\divv\,u)=0.
\end{equation}
Multiply the above by $u$ and integrate to obtain that for $0\le t\le T$,
\begin{align}
\left.\int_{\R^3}{\textstyle\frac{1}{2}}\rho|u|^{2}dx\right |_{0}^{t}dx + \int_{0}^{t}\!\!\!\int_{\R^3} &\rho u(\rho^{-1}\nabla P_\gamma-\tilde\rho_\gamma^{-1}\nabla \tilde P_\gamma)dxds
\notag\\ +& \int_{0}^{t}\!\!\!\int_{\R^3}\left[\mu|\nabla u|^{2} + (\xi-\mu)(\mathrm{div}\,u)^{2}\right]dxds=0,\label{3.6}
\end{align}
where the divergence of a matrix is taken row-wise.
Next we define 
\hfill
$$G(\rho)=\int_{\tilde\rho_\gamma}^{\rho}\!\!\!\int_{\tilde\rho_\gamma}^{r}s^{-1}P'_\gamma(s)dsdr,$$
then using the mass equation, the second term on the left side of $\eqref{3.6}$ can be written as follows
\begin{align*}
\int_{0}^{t}\!\!\!\int_{\R^3}&\rho u(\rho^{-1}\nabla P_\gamma-\tilde\rho_\gamma^{-1}\nabla \tilde P_\gamma)dxds\\
&=\int_0^t\!\!\!\int_{\R^3}\rho u\cdot\nabla\left(\int_{\tilde\rho_\gamma}^{\rho}r^{-1}P'_\gamma(r)dr\right)dxds\\
&=\int_0^t\!\!\!\int_{\R^3}\rho_t\left(\int_{\tilde\rho_\gamma}^{\rho}r^{-1}P'_\gamma(r)dr\right)dxds=\int_0^t\!\!\!\int_{\R^3}G(\rho)_s dxds=\int_{\R^3}G(\rho)dx\Big|_0^t.
\end{align*}
Putting the above into \eqref{3.6}, the result follows.
\end{proof}
\bigskip

Next we derive preliminary $L^2$ bounds for $\nabla u$ and $\dot u$:
\bigskip

\noindent{\bf Lemma 2.4}  \em Assume that the hypotheses and notations of  {\rm Lemma~2.3} are in force. Then for $t\in(0,T]$ and $\gamma=1$,
\begin{align}
\sup_{0 < s\le  t}\sigma\int_{\R^3}|\nabla u|^2 dx + \int_{0}^{t}\!\!\!\int_{\R^3}\sigma\rho|\dot{u}|^2 dxds \le M\left[C_0 + H\right],\label{3.7}
\end{align}
and 
\begin{align}\label{3.10}
\sup_{0<s\le t}\sigma^3\int_{\R^3}|\dot{u}|^2dx +  \int_{0}^{t}\!\!\!\int_{\R^3}\sigma^3|\nabla\dot{u}|^2dxds\le M\left[C_0 + H\right];
\end{align}
while for $\gamma>1$,
\begin{align}\label{3.7*}
\sup_{0 < s\le  t}\sigma\int_{\R^3}|\nabla u|^2 dx + \int_{0}^{t}\!\!\!\int_{\R^3}\sigma\rho|\dot{u}|^2 dxds \le M\left[C_0 + H + D\right],
\end{align}
and 
\begin{align}\label{3.10*}
\sup_{0<s\le t}\sigma^3\int_{\R^3}|\dot{u}|^2dx +  \int_{0}^{t}\!\!\!\int_{\R^3}\sigma^3|\nabla\dot{u}|^2dxds\le M\left[C_0 + H\right].
\end{align}
\rm
\begin{proof}
The proofs are nearly the same as those of \rm{(2.9)} and \rm{(2.12)} in \cite{hoff1}, except that $\tilde\rho_\gamma$ is not necessarily a constant. We prove \eqref{3.7} and \eqref{3.7*} as an example. First, \eqref{3.5} can be rewritten as follows
\begin{align}\label{3.8}
\rho\dot u&+\tilde P_\gamma\nabla(P_\gamma \tilde P_\gamma^{-1}-1)+\rho \tilde P_\gamma^{-1}\nabla \tilde P_\gamma(\rho^{-1}P_\gamma-\tilde\rho_\gamma^{-1}\tilde P_\gamma)-\mu\Delta u-(\xi-\mu)\nabla(\divv\,u)=0.
\end{align}
We multiply \eqref{3.8} by $\sigma\dot u$ and integrate to obtain
\begin{align}\label{3.9}
 \int_{0}^{t}\!\!\!\int_{\R^3}&\sigma\rho|\dot u|^2dxds - \mu\int_0^t\!\!\!\int_{\R^3}\sigma\dot u dxds-(\xi-\mu)\int_0^t\!\!\!\int_{\R^3}\sigma\dot u\nabla(\divv\,u)dxds\notag\\
&-\int_{0}^t\!\!\!\int_{\R^3}\sigma\dot u \tilde P_\gamma\nabla(P_\gamma\tilde P_\gamma^{-1}-1)dxds\notag\\
&\qquad\qquad\qquad-\int_{0}^t\!\!\!\int_{\R^3}\sigma\dot u\rho \tilde P_\gamma^{-1}\nabla \tilde P_\gamma(\rho^{-1}P_\gamma-\tilde\rho_\gamma^{-1}\tilde P_\gamma)dxds=0.
\end{align}
The second and the third term on the left side of \eqref{3.9} can be bounded above by \\$\displaystyle{-\frac{\mu}{2}\sigma\int_{\R^3}|\nabla u|^2dx - \frac{\xi-\mu}{2}\sigma\int_{\R^3}(\divv\,u)^2dx + H}$. The forth term on the left side of \eqref{3.9} can be estimated as follows
\begin{align}
\left|-\int_{0}^t\!\!\!\int_{\R^3}\sigma\dot u \tilde P_\gamma\nabla(P_\gamma\tilde P_\gamma^{-1}-1)dxds\right|\le M\left[\sigma\int_{\R^3}|\nabla u||\rho-\tilde\rho_\gamma|dx+\int_0^{1\wedge t}\!\!\!\int_{\R^3}|\nabla u||\rho-\tilde\rho_\gamma|dxds\right.\notag\\
\left.+\int_{0}^t\!\!\!\int_{\R^3}|\nabla u|^2\right],\notag
\end{align}
where the right side is bounded by $M\left[C_0+A\right]$ by Lemma~2.3. Finally for the fifth term, we notice that for $\gamma=1$, 
$$\rho^{-1}P_\gamma-\tilde\rho_\gamma^{-1}\tilde P_\gamma=L-L=0,$$
so that the term vanishes for the case when $\gamma=1$. On the other hand, for $\gamma>1$, using the definition \eqref{3.3} of $D$, we have
\begin{align*}
-\int_{0}^t\!\!\!\int_{\R^3}&\sigma\dot u\rho \tilde P_\gamma^{-1}\nabla \tilde P_\gamma(\rho^{-1}P_\gamma-\tilde\rho_\gamma^{-1}\tilde P_\gamma)dxds\\
&\le M\int_0^t\!\!\!\int_{\R^3}\sigma|\dot u|\phi(\gamma)|\rho-\tilde\rho_\gamma||\nabla\tilde\rho_\gamma|dxds
\le D^\frac{1}{2}\left[\int_0^t\!\!\!\int_{\R^3}\sigma|\dot u|^2dxds\right]^\frac{1}{2}.
\end{align*}
 Therefore \eqref{3.7} and \eqref{3.7*} follow.
\end{proof}
\medskip

The following auxiliary estimates will be applied to bound the the functional $H$: 
\medskip

\noindent{\bf Lemma 2.5} \em Assume that the hypotheses and notations of  {\rm Lemma~2.3} are in force.
Then for \\$0<t\le 1\wedge T$ and $\gamma\ge1$,
\begin{align}\label{3.11}
\sup_{0\le s\le t}\int_{\R^3}|u|^q dx+\int_0^t\!\!\!\int_{\R^3}\left[|u|^{q-2}|\nabla u|^2 \right.+&\left.|u|^{q-4}|\nabla(|u|^2)|^2\right]dxds\\
&\le M\left[C_0^\frac{p-q}{p-2}N^\frac{q-2}{p-2} + C_0\right]\notag.
\end{align}
\rm
\begin{proof} We multiply \eqref{3.8} by $u|u|^{q-2}$ and integrate to obtain that
\begin{align}\label{3.12}
&\int_{\R^3}|u|^q dx+ \int_0^t\!\!\!\int_{\R^3}|u|^{q-2}|\nabla u|^2dxds + \int_0^t\!\!\!\int_{\R^3}|u|^{q-4}|\nabla(|u|^2)|^2dxds\\
&\le M\left[\int_0^t\!\!\!\int_{\R^3}(P_\gamma-\tilde P_\gamma)\divv(|u|^{q-2}u)dxds+\int_0^t\!\!\!\int_{\R^3}|u|^{q-1}|\nabla\tilde\rho_\gamma||P_\gamma-\tilde P_\gamma|dxds\right.\notag\\
&\qquad\qquad\qquad\qquad-\left.\int_0^t\!\!\!\int_{\R^3}\rho\tilde P_\gamma^{-1}(\rho^{-1}P_\gamma-\tilde\rho_\gamma^{-1}\tilde P_\gamma)|u|^{q-2}u dxds+\int_{\R^3}|u_0|^q dx\right].\notag
\end{align}
The last term on the right side of \eqref{3.12} is bounded by $C_0^\frac{p-q}{p-2}N^\frac{q-2}{p-2}$, and the first and second term can be bounded in terms of $\displaystyle{\sup_{0\le s\le t}\int_{\R^3}|\rho-\tilde\rho_\gamma|^2dx}$. For the remaining term we have
\begin{align*}
-\int_0^t\!\!\!\int_{\R^3}\rho\tilde P_\gamma^{-1}(\rho^{-1}P_\gamma-&\tilde\rho_\gamma^{-1}\tilde P_\gamma)|u|^{q-2}u dxds\\
&\le\int_0^t\!\!\!\int_{\R^3}|\nabla\tilde\rho_\gamma|\phi(\gamma)|\rho-\tilde\rho_\gamma||u|^{q-1}dxds\\
&\le M\left[\sup_{0\le s\le t}\int_{\R^3}|\rho-\tilde\rho_\gamma|^2\right]^\frac{1}{q}\left[\sup_{0\le s\le t}\int_{\R^3}|u|^qdx\right]^\frac{q-1}{q}.
\end{align*}
Using the above on \eqref{3.12} and absorbing terms, the result follows.
\end{proof}
\medskip

We obtain bounds for $u,\omega$ and $ F$ in $W^{1,r}$, these being required for the derivation of estimates for the auxiliary functionals $H$ and $D$:\medskip

\noindent{\bf Lemma 2.6} \em  Let $\gamma\ge1$ and assume that the hypotheses and notations of  {\rm Lemma~2.3} are in force.Then for $r\in(1,\infty)$ and $t\in(0,T]$,
\begin{align}\label{3.13}
||\nabla u(\cdot,t)||_{L^r}\le M\left[|| F(\cdot,t)||_{L^r}+||\omega(\cdot,t)||_{L^r}+||(P_\gamma-\tilde P_\gamma)(\cdot,t)||_{L^r}\right].
\end{align}
The constant $M$ in \eqref{3.13} may depend additionally on $r$.
\rm
\begin{proof} We have from the definition \eqref{1.3} of $F$ that
\begin{align}\label{3.15}
\xi\Delta u^j= F_{x_j}+\xi\omega^{j,k}_{x_k}+(P_\gamma-\tilde P_\gamma)_{x_j}.
\end{align}
Differentiating and taking the Fourier transform we then obtain
\begin{align*}
\xi\widehat{u}^j_{x_l}(y,t)=\frac{y_j y_l}{|y|^2}\widehat{F}(y,t)+\xi\frac{y_k y_l}{|y|^2}\widehat{\omega^{j,k}}(y,t)+\frac{y_k y_l}{|y|^2}(\widehat{P_\gamma-\tilde P_\gamma})(y,t).
\end{align*}
The result \eqref{3.13} then follows immediately from the Marcinkiewicz multiplier theorem (Stein \cite{stein}, pg. 96). 
\end{proof}

\noindent{\bf Lemma 2.7} \em  Assume that the hypotheses and notations of  {\rm Lemma~2.3} are in force. Then for $\gamma>1$,

\begin{align}\label{3.14}
||\nabla F(\cdot,t)||_{L^2}&+||\nabla\omega(\cdot,t)||_{L^2}\notag\\
&\le M\left[||\dot u(\cdot,t)||_{L^2}+||\nabla u(\cdot,t)||_{L^2}+|||\nabla\tilde P||\rho-\tilde\rho_\gamma(\cdot,t)||_{L^2}\right];
\end{align}
while for $\gamma=1$,
\begin{align}\label{3.14*}
||\nabla\left(\frac{F}{\tilde\rho}\right) (\cdot,t)||_{L^2}+||\nabla\left(\frac{\omega}{\tilde\rho}\right)(\cdot,t)||_{L^2}\le M\left[||\dot u(\cdot,t)||_{L^2}+||\nabla u(\cdot,t)||_{L^2}\right].
\end{align}
\rm

\begin{proof}For \eqref{3.14}, using \eqref{1.5} we have
\begin{align}\label{3.16}
\Delta F=(\rho\dot u^j)_{x_j}+[(\tilde P_\gamma)_{x_j}\tilde\rho^{-1}(\tilde\rho_\gamma-\rho)]_{x_j}.
\end{align}
Similar to \eqref{3.15}, we differentiate and take Fourier transform on \eqref{3.16} and apply the multiplier theorem, we get the bounds for $\nabla F$ and similarly for $\divv\,\omega$. To show \eqref{3.14*}, using \eqref{1.3} and the definition \eqref{1.4} of $F$, for $\gamma=1$,
\begin{align}\label{1.3*}
\rho\dot{u}^{j}&=\left(\frac{F}{L\tilde\rho}L\tilde\rho\right)_{x_j}+\mu\omega_k+\rho\tilde\rho^{-1}L\tilde\rho_{x_j}-L\tilde\rho_{x_j}\notag\\
&=\tilde\rho\left(\frac{F}{\tilde\rho}\right)_{x_j}+L\tilde\rho_{x_j}\left[\frac{\xi}{L\tilde\rho}\divv\,u-\rho\tilde\rho^{-1}+1\right]+\mu\left(\frac{\omega}{L\tilde\rho}L\tilde\rho\right)_{x_k}+L\tilde\rho_{x_j}(\rho\tilde\rho^{-1}-1)\notag\\
&=\tilde\rho\left(\frac{F}{\tilde\rho}\right)_{x_j}+\frac{\mu}{L}\left(\frac{\omega}{\tilde\rho}\right)_{x_k}+\frac{\tilde\rho_{x_j}\xi}{\tilde\rho}\divv\,u+\frac{\tilde\rho_{x_j}\mu}{\tilde\rho}\omega,
\end{align}
and so,
\begin{align}\label{3.15*}
||\nabla \left(\frac{F}{\tilde\rho}\right)(\cdot,t)||_{L^2}\le M\left[||\dot u(\cdot,t)||_{L^2}+||\nabla u(\cdot,t)||_{L^2}+||\nabla\left(\frac{\omega}{\tilde\rho}\right)(\cdot,t)||_{L^2}\right].
\end{align}
Using similar method, we obtain
\begin{align}\label{3.16*}
||\nabla\left(\frac{\omega}{\tilde\rho}\right)(\cdot,t)||_{L^2}\le M\left[||\dot u(\cdot,t)||_{L^2}+||\nabla u(\cdot,t)||_{L^2}\right].
\end{align}
Therefore \eqref{3.14*} follows from \eqref{3.15*} and \eqref{3.16*}.
\end{proof}
\medskip

Next we derive a bound for the functional $D$ defined above in \eqref{3.3}:
\medskip

\noindent{\bf Lemma 2.8} \em Assume that the hypotheses and notations of  {\rm Lemma~2.3} are in force. Then for $\gamma=1$ and $t\in(0,T]$,
\begin{align}\label{3.24}
D(t)\le M[C_0+A(t)].
\end{align}
For $\gamma>1$, \eqref{3.24} also holds when $|{\rm supp}(\psi)|$ is sufficiently small.
\rm
\begin{proof}
In view of Lemma~2.3, it suffices to consider $t>1$. We first consider the case when $\gamma=1$. Using the mass equation,
\begin{align}\label{3.21}
\xi\frac{D}{Dt}(\rho-\tilde\rho_\gamma)+\rho(P_\gamma-\tilde P_\gamma)=-\rho\tilde\rho\frac{F}{\tilde\rho} - \xi u\cdot\nabla\tilde\rho_\gamma,
\end{align}
so that we multiply the above by $(\rho-\tilde\rho_\gamma)\rm{sgn}(\rho-\tilde\rho_\gamma)|\nabla\tilde P|^{2}$ and integrate to obtain
\begin{align}\label{3.23}
\int_{\R^3}|\rho-\tilde\rho_\gamma|^2|\nabla\tilde P|^{2}dx+&M^{-1}\int_1^t\!\!\!\int_{\R^3}|\rho-\tilde\rho_\gamma|^2|\nabla\tilde P|^{2}dxds\notag\\
&\le M\left[\int_1^t\!\!\!\int_{\R^3}(|\frac{F}{\tilde\rho}|^2+|u|^2)(1+|x|)^{-2}dxds\right]\notag\\
&\le M\left[\int_1^t\!\!\!\int_{\R^3}(|\nabla \left(\frac{F}{\tilde\rho}\right)|^2+|\nabla u|^2)dxds\right],
\end{align}
where the last line follows by \eqref{1.5.5}. Using \eqref{3.4} and \eqref{3.14*}, the right side of \eqref{3.23} is bounded by $\dis M\left[C_0+A\right]$, and hence \eqref{3.24} holds for $\gamma=1$. 

For $\gamma>1$, using \eqref{1.5.5}, \eqref{3.14} and \eqref{3.23},
\begin{align}\label{3.23*}
\int_1^t\!\!\!\int_{\R^3}|\rho-\tilde\rho_\gamma|^2|\nabla\tilde P|^{2}dxds&\le M\left[C_0+\int_1^t\!\!\!\int_{\R^3}(|F|^2+|u|^2)|\nabla\tilde P|^{2}dxds\right]\notag\\
&\le M\left[C_0+\int_1^t\!\!\!\int_{\R^3}|\nabla u|^2dxds+\int_1^t\!\!\!\int_{\R^3}|F|^2 |\nabla\tilde P|^{2}dxds\right]\notag\\
&\le M\left[C_0+\int_1^t\!\!\!\int_{{\rm supp}(\psi)}|F|^2 |\nabla\tilde P|^{2}dxds\right],
\end{align}
where the last inequality follows from the fact that $|\nabla\tilde P|=\tilde\rho|\nabla\psi|$. Assume that $|{\rm supp}(\psi)|<\infty$, the last term on the right side of \eqref{3.23*} can be estimated as follows
\begin{align*}
\int_1^t\!\!\!\int_{{\rm supp}(\psi)}|F|^2& |\nabla\tilde P|^{2}dxds\le|{\rm supp}(\psi)|^{\frac{1}{3}}\int_1^t\!\!\left[\int_{{\rm supp}(\psi)}|F|^3|\nabla\tilde P|^3 dx\right]^{\frac{2}{3}}ds\notag\\
&\le |{\rm supp}(\psi)|^{\frac{1}{3}}\int_1^t\left[\int_{\R^3}|F|^2|\nabla\tilde P|^2 dx\right]^{\frac{1}{2}}\left[\int_{\R^3}|\nabla F|^2|\nabla\tilde P|^2 dx+\int_{\R^3}|F|^2|D^2_{x}\tilde P|^2\right]^{\frac{1}{2}}ds\notag\\
&\le |{\rm supp}(\psi)|^{\frac{1}{3}}\left[\int_1^t\!\!\!\int_{\R^3}|F|^2|\nabla\tilde P|^2 dxds+M\int_1^t\!\!\!\int_{\R^3}(|\dot u|^2+|\nabla u|^2+|\nabla\tilde P|^2|\rho-\tilde\rho|^2)dxds\right]\notag\\
&\le |{\rm supp}(\psi)|^{\frac{1}{3}}\left[\int_1^t\!\!\!\int_{\R^3}|F|^2|\nabla\tilde P|^2 dxds+M\int_1^t\!\!\!\int_{\R^3}|\nabla\tilde P|^2|\rho-\tilde\rho|^2)dxds+M(C_0+A)\right],
\end{align*}
hence if $|{\rm supp}(\psi)|$ is sufficiently small, we have
\begin{align*}
\int_1^t\!\!\!\int_{\R^3}|\rho-\tilde\rho_\gamma|^2|\nabla\tilde P|^{2}dxds\le M\left[C_0+A\right]
\end{align*}
and \eqref{3.24} follows.
\end{proof}
\medskip

We can now obtain the required estimates for the functionals $H$ defined above in \eqref{3.2} :
\medskip

\noindent{\bf Lemma 2.9}  \em Assume that the hypotheses and notations of  {\rm Lemma~2.3} are in force. Then for $\gamma=1$ and $t\in(0,T]$,
\begin{align}\label{3.17}
H(t)\le M\left[C_0+C_0^{2}+A(t)^2\right].
\end{align}
For $\gamma>1$, \eqref{3.17} also holds when $|{\rm supp}(\psi)|$ is sufficiently small.
\rm
\begin{proof}
We first consider the case for $\gamma=1$. In view of the definition \eqref{3.2} of $H$, we bound only the term \newline $\displaystyle{\int_{0}^{t}\!\!\!\int_{\R^3}\Big[\sigma^{3/2}|\nabla u|^3 + \sigma^3 |\nabla u|^4 \Big]dxds}$. \newline The remaining term $\dis\sum_{1\le k_i,j_m\le 3}\Big|\int_{0}^{t}\!\!\!\int_{\R^3}\sigma u^{j_1}_{x_{k_1}}u^{ j_2}_{x_{k_2}}u^{j_3}_{x_{k_3}}dxds\Big |$ is bounded exactly as in Hoff \cite{hoff1} pp. 29--32.

First we have from Lemma~2.6 that
\begin{align}\label{3.18}
\int_{0}^{t}\!\!\!\int_{\R^3}\sigma^3|\nabla u|^4 dxds\le M\left[\int_{0}^{t}\!\!\!\int_{\R^3}\sigma^3\left(|\rho-\tilde\rho_\gamma|^4+|F|^4+|\omega|^4\right)dxds\right].
\end{align}
Applying \eqref{1.5.1} we can bound the second term on the right by
\begin{align*}
\bar\rho^{4}\int_{0}^{t}\!\!\!\int_{\R^3}\sigma^3|\frac{F}{\tilde\rho}|^4 dxds\le\left(\sup_{0\le s\le t}\int_{\R^3}\sigma|\frac{F}{\tilde\rho}|^2 dx\right)^\frac{1}{2}&\left(\sigma^3\int_{\R^3}|\nabla\left(\frac{F}{\tilde\rho}\right)|^2 dx\right)^\frac{1}{2}\\
&\times\left(\int_{0}^{t}\!\!\!\int_{\R^3}\sigma|\nabla\left(\frac{F}{\tilde\rho}\right)|^2 dxds\right).
\end{align*}
Now from the definition of $ F$ and Lemma~2.3, $$\left(\sup_{0\le s\le t}\int_{\R^3}\sigma| \frac{F}{\tilde\rho}|^2 dx\right)^{\frac{1}{2}}\le M(C_0+A)^\frac{1}{2}.$$ Also, from \eqref{3.14*},
\begin{align*}
\int_0^t\!\!\!\int_{\R^3}\sigma|\nabla \left(\frac{F}{\tilde\rho}\right)|^2dxds&\le M\left[\int_0^t\!\!\!\int_{\R^3}\sigma\left(|\dot{u}|^2+|\nabla u|^2\right)dxds\right]\le M\left[A+C_0\right]
\end{align*}
and 
\begin{align*}
\sigma^3\int_{\R^3}|\nabla F|^2 dx&\le M\left[\sup_{0\le s\le t}\int_{\R^3}\sigma^3\left(|\dot{u}|^2+|\nabla u|^2\right)dx\right]\le MA.
\end{align*}
Thus
\begin{align}\label{3.19}
\int_0^t\!\!\!\int_{\R^3}\sigma^3|F|^4\le M\left[C_0+A\right]^2.
\end{align}
Applying the results of  Lemma~2.7 in a similar way, we obtain that
\begin{align}\label{3.20}
\bar\rho^{4}\int_{0}^{t}\!\!\!\int_{\R^3}\sigma^3|\frac{\omega}{\tilde\rho}|^4 dxds&\le\left(\sup_{0\le s\le t}\int_{\R^3}\sigma|\frac{\omega}{\tilde\rho}|^2 dx\right)^\frac{1}{2}\left(\sigma^3\int_{\R^3}|\nabla\left(\frac{\omega}{\tilde\rho}\right)|^2 dx\right)^\frac{1}{2}\notag\\
&\qquad\qquad\qquad\qquad\qquad\,\,\,\,\,\,\times\left(\int_{0}^{t}\!\!\!\int_{\R^3}\sigma|\nabla\left(\frac{\omega}{\tilde\rho}\right)|^2 dxds\right)\notag\\
&\le M\left[C_0+A\right]^2.
\end{align}
For the first term on the right of \eqref{3.18}, multiply \eqref{3.22} by $\sigma^3|\rho-\tilde\rho_\gamma|^3\rm{sgn}(\rho-\tilde\rho_\gamma)$, integrate and use \eqref{3.19},
\begin{align}\label{3.22}
\int_0^t\!\!\!\int_{\R^3}\sigma^3|\rho-\tilde\rho_\gamma|^4dxds&\le M\left[C_0+\int_0^t\!\!\!\int_{\R^3}\sigma^3(|F^4+|u|^4)dxds\right]\notag\\
&\le M\left[C_0+\left(C_0+A\right)^2\right].
\end{align}
Hence we obtain a bound for $\dis\int_{0}^{t}\!\!\!\int_{\R^3}\sigma^3|\nabla u|^4 dxds$ from \eqref{3.19}, \eqref{3.20} and \eqref{3.22}. Bounds for the term $\dis\int_{0}^{t}\!\!\!\int_{\R^3}\sigma^{3/2}|\nabla u|^3dxds$ are obtained in a similar way, thereby proving \eqref{3.17} for $\gamma=1$.

For $\gamma>1$, using \eqref{3.14},
\begin{align*}
\int_{0}^{t}\!\!\!\int_{\R^3}\sigma^3|F|^4 dxds&\le\left(\sup_{0\le s\le t}\int_{\R^3}\sigma| F|^2 dx\right)^\frac{1}{2}\left[\int_{0}^{t}\left(\int_{\R^3}\sigma|\nabla F|^2 dx\right)^{\frac{3}{2}}ds\right]\\
&\le(C_0+A)^{\frac{1}{2}}\left[\int_{0}^{t}\left(\int_{\R^3}\sigma(|\dot u|^2+|\nabla u|^2+|\rho-\tilde\rho|^2|\nabla\tilde P|^2 dx\right)^{\frac{3}{2}}ds\right]\\
&\le M(C_0+A)^{\frac{1}{2}}(A^{\frac{3}{2}}+A^{\frac{1}{2}}C_0+C_0^{\frac{1}{2}}D)\\
&\le M(C_0+A+D)^2,
\end{align*}
and similarly,
\begin{align*}
\int_0^t\!\!\!\int_{\R^3}\sigma^3|\rho-\tilde\rho_\gamma|^4dxds+\int_{0}^{t}\!\!\!\int_{\R^3}\sigma^3|\omega|^4 dxds\le M(C_0+A+D)^2.
\end{align*}
So when $|{\rm supp}(\psi)|$ is sufficiently small, we can apply \eqref{3.24} to conclude that
\begin{align*}
\int_{0}^{t}\!\!\!\int_{\R^3}\sigma^3|\nabla u|^4 dxds\le M[C_0+(C_0+A)^2]
\end{align*}
which proves \eqref{3.17} for $\gamma>1$.
\end{proof}
\bigskip

\begin{proof}[{\bf Proof of Theorem~2.1 and 2.2:}] Theorem 2.1-2.2 now follows immediately from the bounds \eqref{3.4}, \eqref{3.7}, \eqref{3.10}, \eqref{3.11}, \eqref{3.24} and \eqref{3.17}, and the fact that the functional $A$ is continuous in time.\end{proof}

\section{Pointwise bounds for the density}
\bigskip

In this section we derive pointwise bounds for the density $\rho$ for $\gamma\ge1$, bounds which are independent both of time and of initial smoothness. This will then close the estimates of Theorem~2.1-2.2 to give an uncontingent estimate for the functional $A$ defined in \eqref{3.1}. The result is as follows:
\medskip

\noindent{\bf Theorem~3.1} \em  Let $\gamma=1$. Assume that the system parameters in {\rm \eqref{1.1}} satisfy the conditions in {\rm (1.6) and (1.9)-\eqref{1.1.16}}. Let positive numbers $\mathcal{C}, L,\rho_{\infty},\bar\rho,\underline\rho,\delta$ be given satisfying \eqref{1.6}-\eqref{1.7} and \eqref{1.9.1}. Then given positive numbers $N$ and $d<\delta$, there are positive constants $a, M, \theta$ depending on the parameters and assumptions in {\rm (1.6) and (1.9)-\eqref{1.1.16}}, on $\mathcal{C}, L,\rho_{\infty},\bar\rho,\underline\rho$, on $N$ and on a positive lower bound for $d$, such that if $(\rho-\tilde\rho_\gamma,u)$ is a solution of \eqref{1.1}-\eqref{1.2} on $\R^3\times[0,T]$ with initial data $(\rho_0-\tilde\rho_\gamma,u_0)\in H^3(\R^3))$ satisfying \eqref{1.2.1}-\eqref{1.2.35} with $C_0 <a$, then in fact
\begin{align*}
\mbox{$\underline\rho\le\rho(x,t)\le\bar\rho$ on $\R^3\times[0,T]$},
\end{align*}
and
\begin{align*}
A(T)\le MC_{0}^{\theta}.
\end{align*}
\rm

\noindent{\bf Theorem~3.2} \em  Let $\gamma>1$. Assume that the system parameters in {\rm \eqref{1.1}} satisfy the conditions in {\rm (1.6) and (1.9)-\eqref{1.1.16}}. Let positive numbers $\mathcal{C}, L,\rho_{\infty},\bar\rho,\underline\rho,\delta$ be given satisfying \eqref{1.6}-\eqref{1.7} and \eqref{1.9.2}. Then given positive numbers $N$ and $d<\delta$, there are positive constants $a, M, \theta$ depending on the parameters and assumptions in {\rm (1.6) and (1.9)-\eqref{1.1.16}}, on $\gamma,\mathcal{C}, L,\rho_{\infty},\bar\rho,\underline\rho$, on $N$ and on a positive lower bound for $d$ such that if $|{\rm supp}(\psi)|<a$ and if $(\rho-\tilde\rho_\gamma,u)$ is a solution of \eqref{1.1}-\eqref{1.2} on $\R^3\times[0,T]$ with initial data $(\rho_0-\tilde\rho_\gamma,u_0)\in H^3(\R^3))$ satisfying \eqref{1.2.1}-\eqref{1.2.35} with $C_0 <a$, then in fact
\begin{align*}
\mbox{$\underline\rho\le\rho(x,t)\le\bar\rho$ on $\R^3\times[0,T]$},
\end{align*}
and
\begin{align*}
A(T)\le MC_{0}^{\theta}.
\end{align*}
\rm

\begin{proof} We prove Theorem~3.1 as an example, the proof of Theorem~3.2 is just similar. Most of the details are reminiscent of those in Hoff-Suen \cite{hoffsuen} section 3. First, we choose positive numbers $b$ and $b'$ satisfying
\begin{align*}
\underline\rho<b<\underline\rho+d<\bar\rho-d<b'<\bar\rho.
\end{align*}
Recall that $\rho_0$ takes values in $[\underline\rho+d,\bar\rho-d]$, so that $\rho\in[\underline\rho,\bar\rho]$ on $\R^3\times [0,\tau]$ for some positive $\tau$. It then follows from Theorem~2.1 that
$A(\tau)\le MC_0^\theta$, where $M$ is now fixed.
We shall show that if $C_0$ is further restricted, then in fact $b< \rho< b'$ on $\R^3\times [0,\tau]$, and so by a simple open-closed argument that $b<\rho<b' $ on all of $\R^3\times[0,T]$, we have $A(T)\le MC_0^\theta$ as well. We shall prove the required upper bound, the proof of the lower bound being similar. 

Fix $y\in\R^3$ and define the corresponding particle path $x(t)$ by
\begin{align*}
\left\{ \begin{array}
{lr} \dot{x}(t)
=u(x(t),t)\\ x(0)=y.
\end{array} \right.
\end{align*}
Suppose that there is a time $t_1\le\tau$ such that $\rho(x(t_1),t_1)=b'$. We may take $t_1$ minimal and then choose 
$t_0<t_1$ maximal such that $\rho(x(t_0),t_0) =\bar\rho-d$. Thus $\rho(x(t),t)\in [(\bar\rho-d), b']$ for $t\in [t_0,t_1]$.  We consider two cases:

\medskip

\noindent Case 1: $t_0<t_1\le T\wedge1$.

We have from the definition \eqref{1.3} of $F$ and the mass equation that
\begin{align*}
\mu\frac{d}{dt}[\log\rho(x(t),t)]+P(\rho(x(t),t))-P_\gamma(\tilde\rho_\gamma(x(t))=-F(x(t),t).
\end{align*}
Integrating from $t_0$ to $t_1$ and and abbreviating $\rho(x(t),t)$ by $\rho(t)$, etc., we then obtain
\begin{align}\label{3.1.02}
\mu[\log(b')-\log(\bar\rho-d)]+\int_{ t_0}^{t_1}[P_\gamma(s)-\tilde P_\gamma(s)]ds=-\int_{ t_0}^{t_1}F(s)ds.
\end{align}
We shall show that 
\begin{align}\label{3.1.03}
-\int_{ t_0}^{t_1}F(s)ds\le\tilde MC_0^{\theta}
\end{align}
for a constant $\tilde M$ which depends on the same quantities as the $M$ from Theorem~2.1 (which has been fixed).
If so, then from \eqref{3.1.02}, 
\begin{align}\label{3.1.04}
\mu[\log(b')-\log(\bar\rho-d)]\le\tilde MC_0^{\theta},
\end{align}
where the last inequality holds because $P_\gamma(s)-\tilde P_\gamma(s)$ is nonnegative on $[t_0,t_1]$.  But \eqref{3.1.04} cannot hold if $C_0$ is small depending on $\tilde M, b',$ and $\bar\rho-d$. Stipulating this smallness condition, we therefore conclude that there is no time $t_1$ such that $\rho(t_1) = \rho(x(t_1),t_1) = b'$. Since $y\in \R^3$ was arbitrary, it follows that $\rho<b'$ on $\R^3\times [0,\tau]$, as claimed. The proof that $b<\rho$ is similar.

To prove \eqref{3.1.03} we let $\Gamma$ be the fundamental solution of the Laplace operator in $\R^3$ and apply \eqref{1.5} to write
\begin{align}
-\int_{ t_0}^{t_1}F(s)ds&=-\int_{ t_0}^{t_1}\!\!\!\int_{\R^3}\Gamma_{x_j}(x(s)-y)\rho\dot{u}^{j}(y,s)dyds\notag\\
&\qquad-\int_{ t_0}^{t_1}\!\!\!\int_{\R^3}\Gamma_{x_j}(x(s)-y)\left[(\tilde P_\gamma)_{x_j}\tilde\rho_\gamma^{-1}(\tilde\rho_\gamma-\rho)\right]dyds.\label{3.1.05}
\end{align}
The first integral on the right here is bounded exactly as in  Lemma~4.2 of Hoff~\cite{hoff1}:
\begin{align*}
\int_{ t_0}^{t_1}\!\!\!\int_{\R^3}&\Gamma_{x_j}(x(t)-y)\rho\dot{u}^{j}(y,t)dyds\\
&\le||\Gamma_{x_j}*(\rho u^{j})(\cdot,t_1)||_{L^{\infty}}+||\Gamma_{x_j}*(\rho u^{j})(\cdot, t_0)||_{L^{\infty}}\\
&\qquad+\left |\int_{ t_0}^{t_1}\!\!\!\int_{\R^3}\Gamma_{x_j x_k}(x(s)-y)\left[u^{k}((x(s),s)-u^{k}(y,s)\right](\rho u^{j})(y,s)dyds\right |\\
&\le \tilde MC_0^{\theta},
\end{align*}
and the second integral on the right side of \eqref{3.1.05} can be bounded as follows
\begin{align*}
-\int_{ t_0}^{t_1}\!\!\!\int_{\R^3}\Gamma_{x_j}&(x(s)-y)[(\tilde P_\gamma)_{x_j}\tilde\rho_\gamma^{-1}(\tilde\rho_\gamma-\rho)]dyds\\
&\le\tilde M\int_{0}^{1}\!\!\!\int_{\R^3}|y|^{-2}|\rho-\tilde\rho_\gamma|^2 dyds\\
&\le\tilde M\int_0^1\left[\int_{|y|\le1}|y|^{-\frac{5}{2}}dy\right]^\frac{4}{5}\left[\int_{|y|\le1}|\rho-\tilde\rho_\gamma|^5 dy\right]^\frac{1}{5}ds\\
&\qquad\qquad+\tilde M\int_0^1\left[\int_{|y|>1}|y|^{-4}dy\right]^\frac{1}{2}\left[\int_{|y|>1}|\rho-\tilde\rho_\gamma|^2 dy\right]^\frac{1}{2}ds\\
&\le\tilde MC_0^{\theta},
\end{align*}
where the last inequality follows from Theorem~2.1. Thus \eqref{3.1.03} is proved.

\medskip

\noindent Case 2: $1\le t_0<t_1$.

Again by the mass equation and the definition \eqref{1.3} of $F$,
\begin{align*}
\frac{d}{dt}(\rho(t)-\tilde\rho_\gamma(t))+\mu^{-1}\rho(t)&(P_\gamma(t)-\tilde P_\gamma(t))\\
&=-\mu^{-1}\rho(t)F(t)-\mu^{-1}u(t)\cdot\nabla\tilde\rho_\gamma(t).
\end{align*}
Multiplying by $(\rho(t)-\tilde\rho_\gamma(t))^3$ we get
\begin{align}\label{3.1.07}
&\frac{1}{4}\frac{d}{dt}(\rho(t)-\tilde\rho_\gamma(t))^4 + \mu^{-1}g(t)\rho(t)(\rho(t)-\tilde\rho_\gamma(t))^4\\
& =-\mu^{-1}\rho(t)(\rho(t)-\tilde\rho_\gamma(t))^3F(t)-\mu^{-1}u(t)\cdot\nabla\tilde\rho_\gamma(t)(\rho(t)-\tilde\rho_\gamma(t))^3,\notag
\end{align}
where $g(t)=(P_\gamma(t)-\tilde P_\gamma(t))(\rho(t)-\tilde\rho_\gamma(t))^{-1}\ge0$ on $[t_0,t_1]$. We integrate  \eqref{3.1.07} to obtain
\begin{align}\label{3.1.08}
(\rho(t_1)-\tilde\rho_\gamma(t_1))^4-(\rho(t_0)&-\tilde\rho_\gamma(t_0))^4\\
&\le\tilde M \int_{t_0}^{t_1}\left[||F(\cdot,s)||_{\infty}^4+||u(\cdot,t)||_{\infty}^4\right]ds.\notag
\end{align}
We shall show that
\begin{align}\label{3.1.09}
\tilde M \int_{t_0}^{t_1}\left[||F(\cdot,s)||_{\infty}^4+||u(\cdot,t)||_{\infty}^4\right]ds.\le\tilde MC_0^{\theta},
\end{align}
so that from \eqref{3.1.08},
\begin{align*}
0<(b'-\bar\rho)^4-(\bar\rho-d-\bar\rho)^4\le\tilde MC_0^{\theta}.
\end{align*}
This cannot hold if $C_0$ is sufficiently small, however, so that, as in Case 1, there is no time $t_1$ such that $\rho(t_1) = \rho(x(t_1),t_1) = b'$. Since $y\in \R^3$ was arbitrary, it follows that $\rho<b'$ on $\R^3\times [0,\tau]$, as claimed. 

To prove \eqref{3.1.09}, it suffices to consider the term $\dis\tilde M \int_{t_0}^{t_1}||F(\cdot,s)||_{\infty}^4ds$. We apply \eqref{1.5.4} to get
\begin{align}\label{3.1.10}
\int_{ t_0}^{t_1}||F(\cdot,s)||^{2}_{\infty}ds&\le\int_{ t_0}^{t_1}\left[||\rho\dot u(\cdot,s)||^2_{L^4}+||\nabla\tilde P_\gamma\tilde\rho_\gamma^{-1}(\tilde\rho_\gamma-\rho)(\cdot,s)||^2_{L^4}\right]ds\notag\\
&\qquad+\int_{ t_0}^{t_1}\left[||\rho\dot u(\cdot,s)||^2_{L^2}+||\nabla\tilde P_\gamma\tilde\rho_\gamma^{-1}(\tilde\rho_\gamma-\rho)(\cdot,s)||^2_{L^2}\right]ds,
\end{align}
and the right side of the above is readily seen to be bounded by $\tilde MC_0^{\theta}$. Thus \eqref{3.1.09} is proved.
\end{proof}

\section{Global Existence of Weak Solutions:   Proof of Theorem~1.1-1.2}
\bigskip

In this section we complete the proof of Theorem~1.1 and 1.2 by constructing weak solutions as limits of smooth solutions. Specifically we fix the constants $a$ and $M$ defined in Theorems 2.1-2.2 and 3.1-3.2, we let initial data $(\rho_0,u_0)$ be given  satisfying the hypotheses \eqref{1.2.1}-\eqref{1.2.35} and \eqref{1.3.1} of Theorem~1.1 and Theorem~1.2, and we take $(\rho^\eta_0,u^\eta_0)$ to be smooth approximate initial data obtained by convolving $(\rho_0,u_0)$ with a standard mollifying kernel of width $\eta>0$. We apply the local existence results (see Nash \cite{nash} or Tani \cite{tani}) to show that there is a smooth local solution $(\rho^\eta,u^\eta)$ of \eqref{1.1}-\eqref{1.2} with initial data $(\rho^\eta_0,u^\eta_0)$, defined up to a positive time $T$, which may depend on $\eta$. The {\em a priori} estimates of Theorem~3.1-3.2 then apply to show that
\begin{align}\label{5.0.1}
\mbox{$A(t)\le MC_{0}^{\theta}$ and $\underline{\rho}\le\rho^{\eta}(x,t)\le\bar{\rho}$,}
\end{align}
where $A(t)$ is defined by \eqref{3.1} but with $(\rho,u)$ replaced by $(\rho^\eta,u^\eta)$. By standard arguments together with the bounds \eqref{5.0.1}, $(\rho^\eta,u^\eta)$ exists and satisfies \eqref{1.1}-\eqref{1.2} for all time.

Those bounds in \eqref{5.0.1} will provide the compactness needed to extract the desired solution $(\rho,u)$ in the limit as $\eta\to 0$.
We begin by proving uniform H\"older continuity of the families $\{u^{\eta}\}$ away from $t=0$:
\bigskip

\noindent{\bf Lemma 4.1} \em Given $\tau>0$ there is a constant $C=C(\tau)$ such that, for all $\eta>0$,
\begin{align}\label{5.1.1}
\langle u^\eta(\cdot,t)\rangle^{\frac{1}{2},\frac{1}{4}}_{\R^3\times [\tau,\infty)}\le C(\tau)C_0^\theta.
\end{align}
\rm
\begin{proof} The proof is exactly as in Hoff \cite{hoff1}, pg. 33 and pg. 41--42.
\end{proof}
\bigskip
Compactness of the approximate solutions $(\rho^\eta,u^\eta)$ now follows:
\bigskip

\noindent{\bf Lemma 4.2} \em There is a sequence $\eta_k\to 0$ and functions $u$ and $\rho$ such that as $k\to\infty$,
\begin{align}\label{5.3.1}
\mbox{ $u^{\eta_k}\rightarrow u$ uniformly on compact sets in $\R^3\times (0,\infty)$};\end{align}
\begin{align}\label{5.3.2}
\nabla u^{\eta_k}(\cdot,t),\nabla\omega^{\eta_k}(\cdot,t)\rightharpoonup\nabla u(\cdot,t),\nabla\omega(\cdot,t)
\end{align}
weakly in $L^2(\R^3)$
for all $t>0$; 
\begin{align}\label{5.3.3}
\sigma^{\frac{1}{2}}\dot{u}^{\eta_k},\sigma^{\frac{3}{2}}\nabla\dot{u}^{\eta_k}\rightharpoonup\sigma^{\frac{1}{2}}\dot{u}\sigma^{\frac{3}{2}}\nabla\dot{u}
\end{align}
weakly in $L^2(\R^3\times[0,\infty))$; and
\begin{align}\label{5.3.4}
\rho^{\eta_k}(\cdot,t)\to \rho(\cdot,t)
\end{align}
strongly in $L^2_{loc}(\R^3)$ for every $t\ge 0$.
\rm

\begin{proof} The uniform convergence \eqref{5.3.1} follows from Lemma~4.1 via a diagonal process. The statements in  \eqref{5.3.2} and \eqref{5.3.3} then follow for this same sequence from \eqref{5.0.1} and elementary considerations based on the equality of weak-$L^2$ derivatives and distribution derivatives. The convergence of approximate densities \eqref{5.3.4} for a further subsequence requires a more involved argument, given in Lions \cite{lions} pp. 21--23 and extended by Feireisl \cite{feireisl1}, pp. 63--64 and 118--127.\end{proof}
\bigskip

\begin{proof}[{\bf Proof of Theorem~1.1-1.2:}] We only prove Theorem~1.1 since Theorem~1.2 can be proved in a similar way. It is clear that the limiting functions $(\rho,u)$ of Lemma~4.2 inherit the bounds in \eqref{1.4.2}-\eqref{1.4.8} from \eqref{5.0.1} and \eqref{5.1.1} (but notice that no statements are made in \eqref{1.4.7} concerning $\dot u(\cdot,t)$). It is also clear from the modes of convergence described in Lemma~4.2 that $(\rho,u)$ satisfies the weak forms \eqref{1.2.5}-\eqref{1.2.6} of the differential equations in \eqref{1.1} as well as the initial condition \eqref{1.2}. The continuity statement \eqref{1.4.1} then follows easily from these weak forms together with the bounds in \eqref{1.4.8}. 

It remains to show \eqref{1.30}. By the same argument as in Hoff \cite{hoff1} pp. 44--47, we have for $r_2\in(2,\infty)$,
\begin{align*}
\lim\limits_{t\rightarrow\infty}||\rho(\cdot,t)-\tilde\rho_\gamma(\cdot)||_{L^\infty}+||u(\cdot,t)||_{L^{r_2}}=0.
\end{align*}
So using \eqref{1.4.8}, for $r_1\in(2,\infty)$,
\begin{align*}
\lim\limits_{t\rightarrow\infty}||\rho(\cdot,t)-\tilde\rho_\gamma(\cdot)||_{L^{r_1}}\le CC_0^\theta\lim\limits_{t\rightarrow\infty}||\rho(\cdot,t)-\tilde\rho_\gamma(\cdot)||_{L^\infty}=0.
\end{align*}
This completes the proof of Theorem~1.1.
\end{proof}
\bigskip


\subsection*{Acknowledgment}
I would like to thank David Hoff for his helpful discussions and valuable comments. I am also indebted to the anonymous referee for his effort on my previous draft.


\begin{thebibliography}{1}
\bibitem{danchin} R. Danchin, Global existence in critical spaces for flows of compressible viscous and heat conductive gases, Arch. Ration. Mech. Anal. 160 (2001), no. 1, 1--39.

\bibitem{feireisl1} E. Feireisl, Dynamics of Viscous Compressible Fluids, Oxford Lecture Series in Mathematics and its Applications, 26. Oxford University Press, Oxford, 2004. 

\bibitem{feireisl2} E. Feireisl, Compressible Navier-Stokes equations with a
non-monotone pressure law, J. Diff. Eqns, 184 (2002), pp.
97--108.

\bibitem{hoff1} D. Hoff, Global solutions of the 
Navier-Stokes 
equations 
for
multidimensional, compressible flow with discontinuous 
initial data, J.
Diff. Eqns. 120, no. 1 (1995), 215--254.

\bibitem{hoff2} D. Hoff, Compressible Flow in a Half-Space with Navier Boundary
Conditions, J. Math. Fluid Mech.  7 (2005), 315--338.

\bibitem{hoff4} D. Hoff, Existence of Solutions to a Model for Sparse, One-dimensional  Fluids, to appear in J. Diff. Eqns.

\bibitem{ladyzhenskaya} O.A. Ladyzhenskaya, The mathematical theory of viscous incompressible flow, 2nd ed., Gordon and Breach, New York, 1969.

\bibitem{lions} P.L. Lions, Mathematical Topics in Fluid Mechanics, vol. 2,
Oxford Lecture Series in Mathematics, 10 (1998).

\bibitem{mn1} A. Matsumura and T. Nishida, The initial value
problem for the equations of motion of viscous and heat-conductive
gases, J. Kyoto Univ. 20 (1980), 67--104.

\bibitem{my} A. Matsumura and N. Yamagata, Global weak solutions of the Navier-Stokes equations for multidimensional compressible flow subject to large external potential forces, 
Osaka J. Math. 38 (2001), no. 2, 399--418.

\bibitem{nash} J. Nash, Le Probl\'{e}me de Cauchy pour les \'{e}quations diff\'{e}rentielles dÕun fluide g\'{e}n\'{e}ral, Bull. Soc. Math. France, 90 (1962), 487--497.

\bibitem{stein} E. M. Stein, Singular Integrals and Differentiability Properties of Functions, Princeton Univ. Press, 1970.

\bibitem{hoffsuen} A. Suen and D. Hoff, Global low-energy weak solutions of the equations of 3D compressible magnetohydrodynamics, to appear in Arch. 
Rational Mechanics Ana.

\bibitem{tani} A. Tani, On the first initial-boundary value problem of compressible viscous fluid motion, Publ. Res. Inst. Math. Sci. 13 (1977), 193--253.

\bibitem{ziemer} W. Ziemer, Weakly differentiable functions, Springer-Verlag, 1989.
\end{thebibliography}
\end{document}